\documentclass[amstex,12pt,english,amssymb]{article}

\usepackage{mathtext}
\usepackage[cp1251]{inputenc}
\usepackage[T2A]{fontenc}
\usepackage[english]{babel}
\usepackage[dvips]{graphicx}
\usepackage{amsmath}
\usepackage{amssymb}
\usepackage{amsxtra}
\usepackage{latexsym}
\usepackage{ifthen}

\usepackage[dvips]{graphicx}
\usepackage{epstopdf}
\usepackage{epsfig}
\usepackage{epic}
\usepackage{eepic}

\begin{document}

%\markboth{\centerline{D. KOVTONYUK, I. PETKOV AND V. RYAZANOV}}
%{\centerline{On regular homeomorphisms of the class $W^{1,1}_{\rm
%loc}$ on plane}}

\newcounter{lemma}[section]
\newcommand{\lemma}{\par \refstepcounter{lemma}%
{\bf Lemma \arabic{section}.\arabic{lemma}.}}
\renewcommand{\thelemma}{\thesection.\arabic{lemma}}

\newcounter{corollary}[section]
\newcommand{\corollary}{\par \refstepcounter{corollary}%
{\bf Corollary \arabic{section}.\arabic{corollary}.}}
\renewcommand{\thecorollary}{\thesection.\arabic{corollary}}

\newcounter{remark}[section]
\newcommand{\remark}{\par \refstepcounter{remark}%
{\bf Remark \arabic{section}.\arabic{remark}.}}
\renewcommand{\theremark}{\thesection.\arabic{remark}}

\newcounter{theorem}[section]
\newcommand{\theorem}{\par \refstepcounter{theorem}%
{\bf Theorem \arabic{section}.\arabic{theorem}.}}
\renewcommand{\thetheorem}{\thesection.\arabic{theorem}}

\newcounter{proposition}[section]
\newcommand{\proposition}{\par \refstepcounter{proposition}%
{\bf Proposition \arabic{section}.\arabic{proposition}.}}
\renewcommand{\theproposition}{\thesection.\arabic{proposition}}

\renewcommand{\theequation}{\arabic{section}.\arabic{equation}}

\def\kohta #1 #2\par{\par\noindent\rlap{#1)}\hskip30pt
\hangindent30pt #2\par}
\def\A{{{\cal {A}}}}
\def\C{{{\Bbb C}}}
\def\Rn{{{\Bbb R}^n}}
\def\lC{{\overline {{\Bbb C}}}}
\def\lRn{{\overline {{\Bbb R}^n}}}
\def\lRm{{\overline {{\Bbb R}^m}}}
\def\lRk{{\overline {{\Bbb R}^k}}}
\def\lBn{{\overline {{\Bbb B}^n}}}
\def\Bn{{{\Bbb B}^n}}
\def\B{{{\Bbb B}}}
\let\text=\mbox

\def\Xint#1{\mathchoice
   {\XXint\displaystyle\textstyle{#1}}%
   {\XXint\textstyle\scriptstyle{#1}}%
   {\XXint\scriptstyle\scriptscriptstyle{#1}}%
   {\XXint\scriptscriptstyle\scriptscriptstyle{#1}}%
   \!\int}
\def\XXint#1#2#3{{\setbox0=\hbox{$#1{#2#3}{\int}$}
     \vcenter{\hbox{$#2#3$}}\kern-.5\wd0}}
\def\dashint{\Xint-}

\def\cc{\setcounter{equation}{0}
\setcounter{figure}{0}\setcounter{table}{0}}

\overfullrule=0pt

\title{{\bf On homeomorphisms with finite distortion \\ in the plane}}

\author{{\bf Denis Kovtonyuk, Igor Petkov  and Vladimir Ryazanov}\\}
\date{\today \hskip 4mm ({\tt KPR121110.tex})}
\maketitle

\large

\begin{abstract} It is shown that every homeomorphism $f$
of finite distortion in the plane is the so-called lower
$Q$-homeomorphism with $Q(z)=K_f(z)$, and, on this base, it is
developed the theory of the boundary behavior of such
homeomorphisms. \end{abstract}

\section{Introduction}

The concept of the generalized derivative was introduced by Sobolev
in \cite{So}. Given a domain $D$ in the complex plane $\C$, the {\bf
Sobolev class} $W^{1,1}(D)$ consists of all functions $f:D\to\C$ in
$L^1(D)$ with first partial generalized derivatives which are
integrable in $D$. A function $f:D\to\C$ belongs to
$W^{1,1}_{\mathrm{loc}}(D)$ if $f\in W^{1,1}(D_*)$ for every open
set $D_*$ with its compact closure $\overline{D_*} \subset D$.

Recall that a homeomorphism $f$ between domains $D$ and $D'$ in
$\C$ is called of {\bf finite distortion} if $f\in
W^{1,1}_{\mathrm{loc}}$ and
\begin{equation}\label{eq1.0KR}||f'(z)||^2\leqslant K(z)\cdot
J_{f}(z)\end{equation} with a.e. finite function $K$ where
$||f'(z)||$ denotes the matrix norm of the Jacobian matrix $f'$ of
$f$ at $z\in D$ and $J_{f}(z)=\det f'(z)$, see \cite{IM}. Later
on, we use the notion $K_{f}(z)$ for the minimal function
$K(z)\geqslant1$ in (\ref{eq1.0KR}). Note that
$||f'(z)||=|f_z|+|f_{\bar{z}}|$ and
$J_f(z)=|f_z|^2-|f_{\bar{z}}|^2$ at the points of total
differentiability of $f$. Thus,
$K_{f}(z)=||f'(z)||^2/J_{f}(z)=\left(|f_z|+|f_{\bar{z}}|\right)/\left(|f_z|-|f_{\bar{z}}|\right)$
if $J_{f}(z)\neq0$, $K_{f}(z)=1$ if $f'(z)=0$, i.e.
$|f_z|=|f_{\bar{z}}|=0$, and $K_{f}(z)=\infty$ at the rest points.

A continuous mapping $\gamma$ of an open subset $\Delta$ of the
real axis $\mathbb{R}$ or a circle into $D$ is called a {\bf
dashed line}, see e.g. Section 6.3 in \cite{MRSY}. Recall that
every open set $\Delta$ in $\mathbb{R}$ consists of a countable
collection of mutually disjoint intervals. This is the motivation
for the term.

Given a family $\Gamma$ of dashed lines $\gamma$ in complex plane
$\C$, a Borel function $\varrho:\C\to[0,\infty]$ is called {\bf
admissible} for $\Gamma$, write $\varrho\in\mathrm{adm}\,\Gamma$,
if
\begin{equation}\label{eq1.2KR}\int\limits_{\gamma}\varrho\,ds\ \geqslant\ 1\end{equation}
for every $\gamma\in\Gamma$. The {\bf (conformal) modulus} of
$\Gamma$ is the quantity \begin{equation}\label{eq1.3KR}M(\Gamma)\
=\
\inf_{\varrho\in\mathrm{adm}\,\Gamma}\int\limits_{\C}\varrho^2(z)\,dm(z)\end{equation}
where $dm(z)$ corresponds to the Lebesgue measure in $\C$. We say
that a property $P$ holds for {\bf a.e.} (almost every)
$\gamma\in\Gamma$ if a subfamily of all lines in $\Gamma$ for
which $P$ fails has the modulus zero, cf. \cite{Fu}. Later on, we
also say that a Lebesgue measurable function
$\varrho:\C\to[0,\infty]$ is {\bf extensively admissible} for
$\Gamma$, write $\varrho\in\mathrm{ext\,adm}\,\Gamma$, if
(\ref{eq1.2KR}) holds for a.e. $\gamma\in\Gamma$, see e.g. 9.2 in
\cite{MRSY}.

The following concept was motivated by Gehring's ring definition
of qua\-si\-con\-for\-ma\-li\-ty in \cite{Ge$_1$}. Given domains
$D$ and $D'$ in $\lC=\C\cup\{\infty\}$,
$z_0\in\overline{D}\setminus\{\infty\}$, and a measurable function
$Q:D\to(0,\infty)$, we say that a homeomorphism $f:D\to D'$ is a
{\bf lower Q-homeomorphism at the point} $z_0$ if
\begin{equation}\label{eq1.4KR}M(f\Sigma_{\varepsilon})\ \geqslant\ \inf\limits_{\varrho\in\mathrm{ext\,adm}\,\Sigma_{\varepsilon}}
\int\limits_{D\cap
R_{\varepsilon}}\frac{\varrho^2(x)}{Q(x)}\,dm(x)\end{equation} for
every ring
$$R_{\varepsilon}=\{z\in\C:\varepsilon<|z-z_0|<\varepsilon_0\},\quad\varepsilon\in(0,\varepsilon_0),\
\varepsilon_0\in(0,d_0)\,,$$ where $$d_0=\sup\limits_{z\in
D}\,|z-z_0|\,,$$ and $\Sigma_{\varepsilon}$ denotes the family of
all intersections of the circles
$$S(r)=S(z_0,r)=\{z\in\C:|z-z_0|=r\},\quad r\in(\varepsilon,\varepsilon_0)\,,$$
with the domain $D$.

The notion can be extended to the case $z_0=\infty\in\overline{D}$
in the standard way by applying the inversion $T$ with respect to
the unit circle in $\lC$, $T(x)=z/|z|^2$, $T(\infty)=0$,
$T(0)=\infty$. Namely, a homeomorphism $f:D\to D'$ is a {\bf lower
$Q$-homeomorphism at} $\infty\in\overline{D}$ if $F=f\circ T$ is a
lower Q$_*$-homeomorphism with $Q_*=Q\circ T$ at $0$. We also say
that a homeomorphism $f:D\to{\lC}$ is a {\bf lower
$Q$-homeomorphism in} $\partial D$ if $f$ is a lower
$Q$-homeomorphism at every point $z_{0}\in\partial D$.

Further we show that every homeomorphism of finite distortion in
the plane is a lower $Q$-homeomorphism with $Q(z)=K_{f}(z)$ and,
thus, the whole theory of the boundary behavior in \cite{KR$_2$},
see also Chapter 9 in \cite{MRSY}, can be applied.

\cc
\section{Preliminaries}

Recall first of all the following topological notion. A domain
$D\subset\C$ is said to be {\bf locally connected at a point}
$z_0\in\partial D$ if, for every neighborhood $U$ of the point
$z_0$, there is a neighborhood $V\subseteq U$ of $z_0$ such that
$V\cap D$ is connected. Note that every Jordan domain $D$ in $\C$
is locally connected at each point of $\partial D$, see e.g.
\cite{Wi}, p. 66.

%%%%%%%%%%%%%%%%%%%% EPS picture %%%%%%%%%%%%%%%%%%%%%%%%%%%%%%%%%
\begin{figure}[h]
\centerline{\includegraphics[scale=0.37]{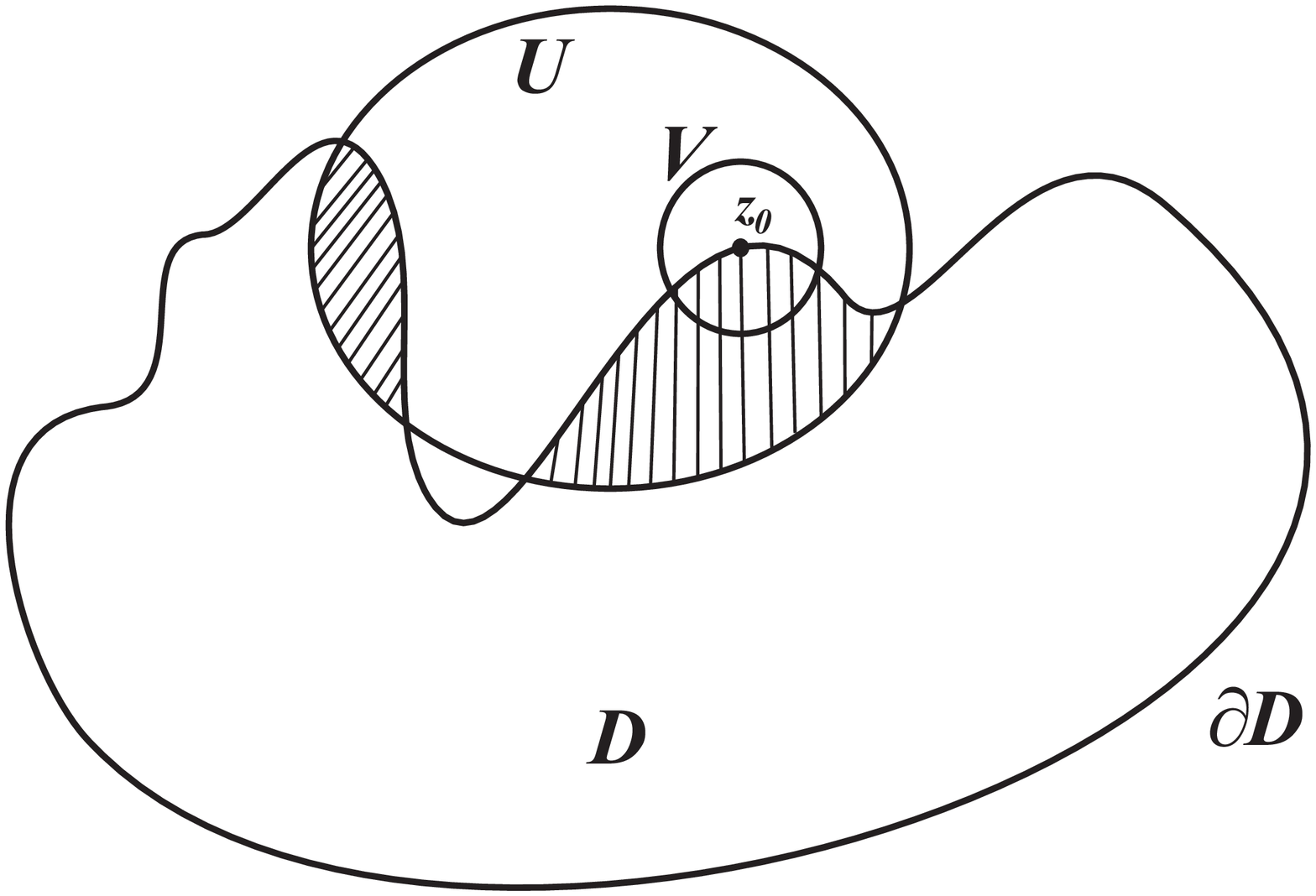}}%
%\centerline{Figure 1.}
\end{figure}
%%%%%%%%%%%%%%%%%%%%%%%%%%%%%%%%%%%%%%%%%%%%%%%%%%%%%%%%%%%%%%%%%%

We say that $\partial D$ is {\bf weakly flat at a point}
$z_0\in\partial D$ if, for every neighborhood $U$ of the point
$z_0$ and every number $P>0$, there is a neighborhood $V\subset U$
of $z_0$ such that
\begin{equation}\label{eq1.5KR}M(\Delta(E,F;D))\geqslant P\end{equation} for
all continua $E$ and $F$ in $D$ intersecting $\partial U$ and
$\partial V$. Here and later on, $\Delta(E,F;D)$ denotes the
family of all paths $\gamma:[a,b]\to\lC$ connecting $E$ and $F$ in
$D$, i.e. $\gamma(a)\in E$, $\gamma(b)\in F$ and $\gamma(t)\in D$
for all $t\in(a,b)$. We say that the boundary $\partial D$ is {\bf
weakly flat} if it is weakly flat at every point in $\partial D$.

%%%%%%%%%%%%%%%%%%%% EPS picture %%%%%%%%%%%%%%%%%%%%%%%%%%%%%%%%%
\begin{figure}[h]
\centerline{\includegraphics[scale=0.4]{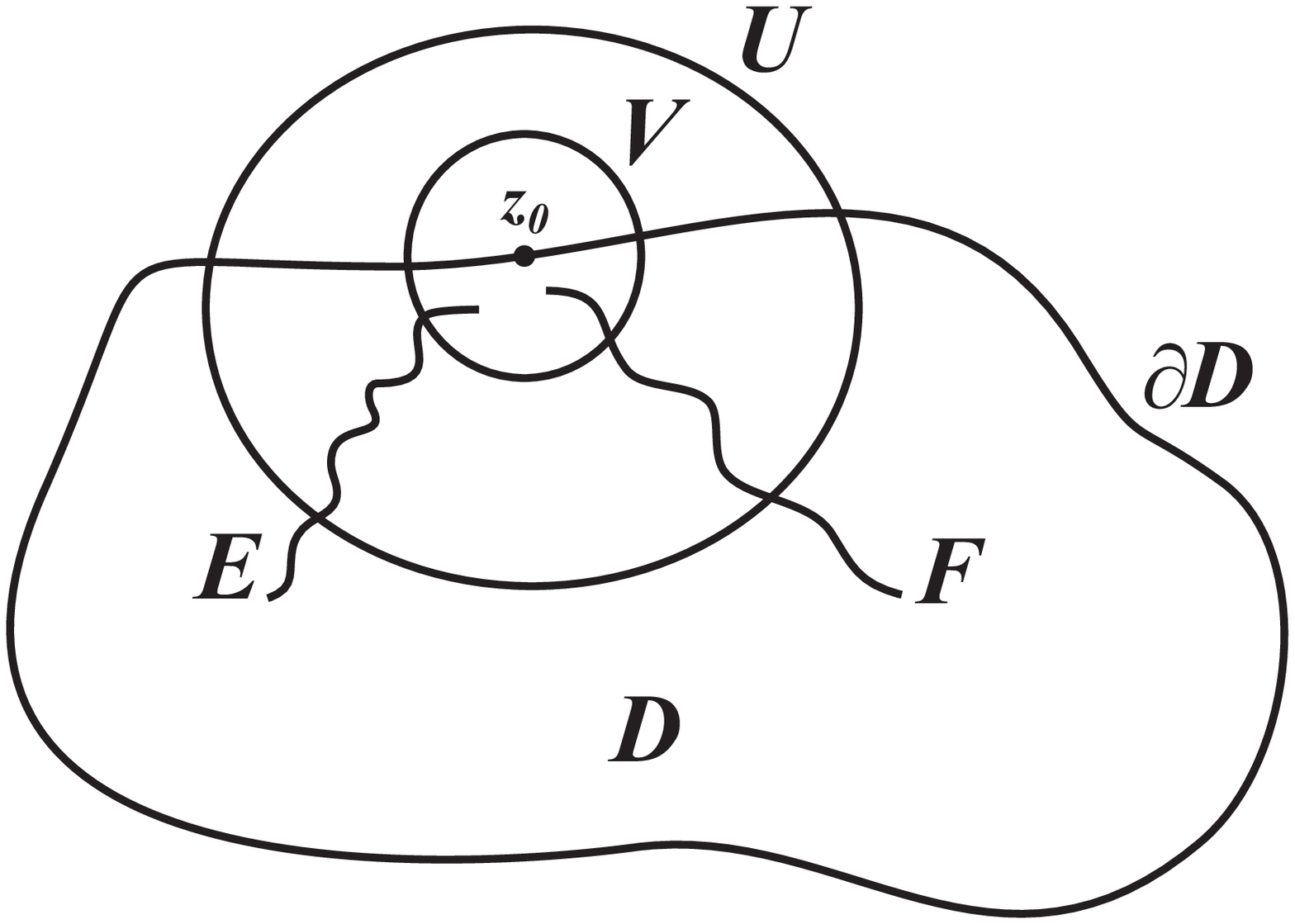}}%
%\centerline{Figure 2.}
\end{figure}
%%%%%%%%%%%%%%%%%%%%%%%%%%%%%%%%%%%%%%%%%%%%%%%%%%%%%%%%%%%%%%%%%%

We also say that a point $z_0\in\partial D$ is {\bf strongly
accessible} if, for every neighborhood $U$ of the point $z_0$,
there exist a compactum $E$ in $D$, a neighborhood $V\subset U$ of
$z_0$ and a number $\delta>0$ such that
\begin{equation}\label{eq1.6KR}M(\Delta(E,F;D))\geqslant\delta\end{equation} for all
continua $F$ in $D$ intersecting $\partial U$ and $\partial V$. We
say that the boundary $\partial D$ is {\bf strongly accessible} if
every point $z_0\in\partial D$ is strongly accessible.

Here, in the definitions of strongly accessible and weakly flat
boundaries, one can take as neighborhoods $U$ and $V$ of a point
$z_0$ only balls (closed or open) centered at $z_0$ or only
neighborhoods of $z_0$ in another fundamental system of
neighborhoods of $z_0$. These conceptions can also be extended in
a natural way to the case of $\lC$ and $z_0=\infty$. Then we must
use the corresponding neighborhoods of $\infty$.

It is easy to see that if a domain $D$ in $\C$ is weakly flat at a
point $z_0\in\partial D$, then the point $z_0$ is strongly
accessible from $D$. Moreover, it was proved by us that if a
domain $D$ in $\C$ is weakly flat at a point $z_0\in\partial D$,
then $D$ is locally connected at $z_0$, see e.g. Lemma 5.1 in
\cite{KR$_2$} or Lemma 3.15 in \cite{MRSY}.

The notions of strong accessibility and weak flatness at boundary
points of a domain in $\C$ defined in \cite{KR$_0$} are
localizations and generalizations of the cor\-res\-pon\-ding
notions introduced in \cite{MRSY$_5$}--\cite{MRSY$_6$}, cf. with
the properties $P_1$ and $P_2$ by V\"ais\"al\"a in \cite{Va$_1$}
and also with the quasiconformal accessibility and the
quasiconformal flatness by N\"akki in \cite{Na$_1$}. Many theorems
on a homeomorphic extension to the boundary of quasiconformal
mappings and their generalizations are valid under the condition
of weak flatness of boundaries. The condition of strong
accessibility plays a similar role for a continuous extension of
the mappings to the boundary. In particular, recently we have
proved the following significant statements, see either Theorem
10.1 (Lemma 6.1) in \cite{KR$_2$} or Theorem 9.8 (Lemma 9.4) in
\cite{MRSY}.

\begin{proposition}{}\label{prKR2.1} {\it Let $D$ and $D'$ be bounded domains in $\C$,
$Q:D\to(0,\infty)$ a measurable function and $f:D\to D'$ a lower
$Q$-homeomorphism in $\partial D$. Suppose that the domain $D$ is
locally connected on $\partial D$ and that the domain $D'$ has a
(strongly accessible) weakly flat boundary. If
\begin{equation}\label{eqKPR2.1}\int\limits_{0}^{\delta(z_0)} \frac{dr}{||\,Q||_{1}(z_0,r)}\
=\ \infty\qquad\forall\ z_0\in\partial D\end{equation} for some
$\delta(z_0)\in(0,d(z_0))$ where $d(z_0)=\sup\limits_{z\in
D}\,|\,z-z_0|$ and
$$||\,Q||_1(z_0,r)=\int\limits_{D\cap
S(z_0,r)}Q(z)\,ds\,,$$ then $f$ has a (continuous) homeomorphic
extension $\overline{f}$ to $\overline{D}$ that maps
$\overline{D}$ (into) onto $\overline{D'}$.}\end{proposition}

Here as usual $S(z_0,r)$ denotes the circle $|z-z_0|=r$.

A domain $D\subset\C$ is called a {\bf quasiextremal distance
domain}, abbr. {\bf QED-domain}, see \cite{GM}, if
\begin{equation}\label{e:7.1}M(\Delta(E,F;\lC)\leqslant K\cdot
M(\Delta(E,F;D))\end{equation} for some $K\geqslant1$ and all
pairs of nonintersecting continua $E$ and $F$ in $D$.

It is well known, see e.g. Theorem 10.12 in \cite{Va$_1$}, that
\begin{equation}\label{eqKPR2.2}M(\Delta(E,F;\C))\geqslant\frac{2}{\pi}\log{\frac{R}{r}}\end{equation}
for any sets $E$ and $F$ in $\C$ intersecting all the circles
$S(z_0,\rho)$, $\rho\in(r,R)$. Hence a QED-domain has a weakly
flat boundary. One example in \cite{MRSY}, Section 3.8, shows that
the inverse conclusion is not true even among simply connected
plane domains.

A domain $D\subset\C$ is called a {\bf uniform domain} if each
pair of points $z_1$ and $z_2\in D$ can be joined with a
rectifiable curve $\gamma$ in $D$ such that
\begin{equation}\label{e:7.2}s(\gamma)\ \leqslant\ a\cdot|\,z_1-z_2|\end{equation}
and \begin{equation}\label{e:7.3}\min\limits_{i=1,2}\
s(\gamma(z_i,z))\ \leqslant\ b\cdot d(z,\partial D) \end{equation}
for all $z\in\gamma$ where $\gamma(z_i,z)$ is the portion of
$\gamma$ bounded by $z_i$ and $z$, see~\cite{MaSa}. It is known that
every uniform domain is a QED-domain but there exist QED-domains
that are not uniform, see~\cite{GM}. Bounded convex domains and
bounded domains with smooth boundaries are simple examples of
uniform domains and, consequently, QED-domains as well as domains
with weakly flat boundaries.

A closed set $X\subset\C$ is called a {\bf null-set for extremal
distances}, abbr. {\bf NED-set}, if
\begin{equation}\label{e:8.1} M(\Delta(E,F;\C))=M(\Delta(E,F;\C\backslash X))\end{equation}
for any two nonintersecting continua $E$ and $F\subset\C\backslash
X$.

\begin{remark}\label{rmKR2.0} It is known that if $X\subset\C$ is a NED-set, then
\begin{equation}\label{e:8.2}|\,X|=0\end{equation} and $X$ does not locally separate $\C$,
see \cite{Va$_2$}, i.e., \begin{equation}\label{e:8.3}
\dim\,X\leqslant 0\,,\end{equation} and hence they are totally
disconnected, see e.g. p. 22 and 104 in \cite{HW}. Conversely, if
a set $X\subset\C$ is closed and is of length zero,
\begin{equation}\label{e:8.4} H^{1}(X)=0\,,\end{equation} then $X$
is a NED-set, see \cite{Va$_2$}. Note also that the complement of
a NED-set in $\C$ is a very particular case of a QED-domain.
\end{remark}

Here $H^{1}(X)$ denotes the 1-dimensional Hausdorff measure
(length) of a set $X$ in $\C$. Also we denote  by $C(X,f)$ the
{\bf cluster set} of the mapping $f:D\to\lC$ for a set
$X\subset\overline D$,
\begin{equation}\label{e:8.5}
C(X,f)\colon=\left\{w\in\lC:w=\lim\limits_{k\to\infty}f(z_k),\
z_k\to z_0\in X,\ z_k\in D\right\}.\end{equation} Note that the
inclusion $C(\partial D,f)\subseteq\partial D'$ holds for every
homeomorphism $f:D\to D'$, see e.g. Proposition 13.5 in \cite{MRSY}.

\cc
\section{The main lemma}

\begin{theorem}\label{thKPR3.1} {\it Let $f:D\to\C$ be a homeomorphism
with finite distortion. Then $f$ is a lower $Q$-homeomorphism at
each point $z_0\in\overline{D}$ with $Q(z)=K_{f}(z)$.}
\end{theorem}

{\it Proof.} Let $B$ be a (Borel) set of all points $z$ in $D$
where $f$ has a total differential with $J_{f}(z)\neq0$ a.e. It is
known that $B$ is the union of a countable collection of Borel
sets $B_l$, $l=1,2,\ldots\,$, such that $f_l=f|_{B_l}$ is a
bi-Lipschitz homeomorphism, see e.g. Lemma 3.2.2 in \cite{Fe}.
With no loss of generality, we may assume that the $B_l$ are
mutually disjoint. Denote also by $B_{*}$ the set of all points
$z\in D$ where $f$ has a total differential with $f'(z)=0$.

Note that the set $B_0=D\setminus(B\cup B_{*})$ has the Lebesgue
measure zero in $\C$ by Gehring--Lehto--Menchoff theorem, see
\cite{GL} and \cite{Me}. Hence by Theorem 2.11 in \cite{KR$_7$}, see
also Lemma 9.1 in \cite{MRSY}, length$(\gamma\cap B_0)=0$ for a.e.
paths $\gamma$ in $D$. Let us show that length$(f(\gamma)\cap
f(B_0))=0$ for a.e. circle $\gamma$ centered at $z_0$.

The latter follows from absolute continuity of $f$ on closed subarcs
of $\gamma\cap D$ for a.e. such circle $\gamma$. Indeed, the class
$W^{1,1}_{\rm loc}$ is invariant with respect to local
quasi-isometries, see e.g. Theorem 1.1.7 in \cite{Maz}, and the
functions in $W^{1,1}_{\rm loc}$ is absolutely continuous on lines,
see e.g. Theorem 1.1.3 in \cite{Maz}. Applying say the
transformation of coordinates $\log(z-z_0)$, we come to the absolute
continuity on a.e. such circle $\gamma$.

Thus, length$(\gamma_{*}\cap f(B_0))=0$ where $\gamma_{*}=f(\gamma)$
for a.e. circle $\gamma$ centered at $z_0$. Now, let
$\varrho_*\in{\mathrm adm}\,f(\Gamma)$ where $\Gamma$ is the
collection of all dashed lines $\gamma\cap D$ for such circles
$\gamma$ and $\varrho_*\equiv0$ outside $f(D)$. Set $\varrho\equiv0$
outside $D$ and
$$\varrho(z)\ \colon=\ \varrho_*(f(z))\left(|\,f_{z}|+|\,f_{\bar{z}}|\right)\quad\ \ {\rm for\ a.e.\ }\ z\in D$$

Arguing piecewise on $B_l$, we have by Theorem 3.2.5 under $m=1$ in
\cite{Fe} that $$\int\limits_{\gamma}\varrho\,ds\ \geqslant\
\int\limits_{\gamma_{*}}\varrho_{*}\,ds_{*}\ \geqslant\ 1\qquad {\rm
for\ a.e.}\ \ \gamma\in\Gamma$$ because length$(f(\gamma)\cap
f(B_0))=0$ and length$(f(\gamma)\cap f(B_{*}))=0$ for a.e.
$\gamma\in\Gamma$, consequently,
$\varrho\in{\mathrm{ext\,adm}}\,\Gamma$.

On the other hand, again arguing piecewise on $B_l$, we have the
inequality
$$\int\limits_{D}\frac{\varrho^2(x)}{K_{f}(z)}\,dm(z)\ \leqslant\
\int\limits_{f(D)}\varrho^2_*(w)\, dm(w)$$ because
$J_{f}(z)=|\,f_{z}|^2-|\,f_{\bar{z}}|^2$ and
$K_{f}(z)=(|\,f_z|+|\,f_{\bar{z}}|)/(|\,f_z|-|\,f_{\bar{z}}|)$ on
$B$ and $K_{f}(z)=1$ and $\varrho(z)=0$ on $B_{*}$. Consequently,
we obtain that
$$M(f\Gamma)\ \geqslant\ \inf\limits_{\varrho\in
{\mathrm{ext\,adm}}\,\Gamma}\int\limits_{D}\frac{\varrho^2(z)}{K_{f}(z)}\,dm(z)\,,$$
i.e. $f$ is really a lower $Q$-homeomorphism with $Q(z)=K_{f}(z)$.

\cc
\section{On the removability of isolated singularities}

In view of Theorem \ref{thKPR3.1} we obtain by Theorem 4.1 in
\cite{KR$_2$} or Theorem 9.3 in \cite{MRSY} the following statement.

\begin{theorem}\label{thKPR4.1} {\it Let $D$ be a domain in $\C$, $z_0\in{D}$,
and $f$ be a homeomorphism with finite distortion of
$D\setminus\{z_0\}$ into $\lC$. Suppose that
\begin{equation}\label{eq8.6.2}\int\limits_{0}^{\varepsilon_0}\frac{dr}{r\cdot
k_f(r)}\ =\ \infty\end{equation} where
$\varepsilon_0<\mathrm{dist}(z_0,\partial D)$ and
\begin{equation}\label{eq8.6.3} k_f(r)\ =\
\dashint_{|\,z-z_0|=r}K_f(z)\ |dz|\ .\end{equation} Then $f$ has a
continuous extension to $D$ in $\lC$.} \end{theorem}

From here we have, in particular, the following consequences.

\begin{corollary}\label{corKPR4.1} {\it Let $D$ be a domain in $\C$
and let $f$ be a homeomorphism with finite distortion of
$D\setminus\{z_0\}$ into $\lC$. If
\begin{equation}\label{eq8.6.5}\dashint_{|\,z-z_0|=r}K_f(z)\ |dz|\ =\
O\left(\log\frac{1}{r}\right)\ \ \ \ \ \ \ \mbox{as}\ r\to 0\
,\end{equation} then $f$ has a continuous extension to $D$ in
$\lC$.}
\end{corollary}

\begin{corollary}\label{corKPR4.2} {\it Let $D$ be a domain in $\C$,
$x_0\in D$, and $f$ be a homeomorphism with finite distortion of
$D\setminus\{z_0\}$ into $\lC$. If
\begin{equation}\label{eq8.6.7}\ \ \ \ \
\dashint_{|\,z-z_0|=r}K_f(z)\ |dz|\ =\
O\left(\log\frac{1}{r}\cdot\log\log\frac{1}{r}\cdot\ldots\cdot\log\ldots\log\frac{1}{r}
\right) \ \ \ \ \ \ \mbox{as}\ r\to0\ , \end{equation}  then $f$ has
a continuous extension to $D$ in $\lC$.} \end{corollary}

\cc
\section{On a continuous extension to boundary points}

In view of Theorem \ref{thKPR3.1} we have by  Theorem 6.1 in
\cite{KR$_2$} or Lemma 9.4 in \cite{MRSY} the next statement.

\begin{lemma}\label{lem5.1} {\it Let $D$ and $D'$ be domains in $\C$,
$z_0\in\partial D$, and $f:D\to D'$ be a homeomorphism with finite
distortion. Suppose that the domain $D$ is locally connected at
$z_0\in\partial D$ and $\partial D'$ is strongly accessible at least
at one point of the cluster set $C(z_0,f)$. If
\begin{equation}\label{eq8.7.4}\int\limits_{0}^{\varepsilon_0}
\frac{dr}{||\,K_f||_{1}(r)}\ =\ \infty\end{equation} where $0\ <\
\varepsilon_0\ <\ d_0\ =\ \sup\limits_{z\in D}\,|\,z-z_0|,$ and
\begin{equation}\label{eq8.7.6} ||\,K_f||_{1}(r)\ =\int\limits_{D\cap
S({z_0},r)}K_f\ ds\ ,\end{equation} then $f$ extends by continuity
to $z_0$ in $\lC$.}\end{lemma}

In particular, we have the following consequence of Lemma
\ref{lem5.1}.

\begin{corollary}\label{thKPR6.1} {\it Let $D$ and $D'$ be {\rm QED} domains in $\C$,
$z_0\in\partial D$, and $f:D\to D'$ be a homeomorphism of finite
distortion. If (\ref{eq8.7.4}) holds, then $f$ extends by continuity
to $z_0$ in $\lC$.} \end{corollary}

Note that the complements of NED sets in $\C$ give very particular
cases of QED domains. Thus, arguing locally, by Theorem
\ref{thKPR6.1}, we obtain the following statement.\newpage

\begin{theorem}\label{thKPR7.1} {\it Let $D$ be a domain in $\C$,
$X\subset D$, and $f$ be a homeomorphism with finite distortion of
$D\backslash X$ into $\lC$. Suppose that $X$ and $C(X,f)$ are {\rm
NED} sets. If
\begin{equation}\label{eq8.9.8} \int\limits_{0}^{\varepsilon_0}
\frac{dr}{||\,K_f||_{1}(r)}\ =\ \infty\end{equation} where
\begin{equation}\label{eq8.9.9} 0\ <\ \varepsilon_0\ <\ d_0\ =\
\mbox{dist}\ (z_0,\partial D) \end{equation} and
\begin{equation}\label{eq8.9.10}
||\,K_f||_{1}(r)\ =\ \int\limits_{|z-z_0|=r}K_f(z)\ |dz|\
,\end{equation} then $f$ can be extended by continuity in $\lC$ to
$z_0$.}
\end{theorem}

\cc
\section{The extension of the inverse mappings to the boundary}

The base of the proof for extending the inverse mappings for
homeomorphisms of finite distortion is the following lemma on the
cluster sets.

\begin{lemma}\label{lemKPR8.1} {\it Let $D$ and $D'$ be domains in $\C$,
$z_1$ and $z_2$ be distinct points in $\partial D$, $z_1\neq\infty$,
and let $f$ be a homeomorphism with finite distortion of $D$ onto
$D'$. Suppose that the function $K_f$ is integrable on the dashed
lines
\begin{equation}\label{eqKPR6.1}D(r)\ =\ \{z\in D:|\,z-z_1|=r\}\
=\ D\cap S(z_1,r)\end{equation} for some set $E$ of numbers
$r<|z_1-z_2|$ of a positive linear measure. If $D$ is locally
connected at $z_1$ and $z_2$ and $\partial D'$ is weakly flat,
then}\end{lemma}
\begin{equation}\label{eq8.10.2} C(z_1,f)\cap C(z_2,f)\ =\
\varnothing. \end{equation}

The of Lemma \ref{lemKPR8.1} follows by Theorem \ref{thKPR3.1} from
Lemma 9.1 in \cite{KR$_2$} or Lemma 9.5 in \cite{MRSY}.

As an immediate consequence of Lemma \ref{lemKPR8.1}, we have the
following statement.

\begin{theorem}\label{thKPR8.2} {\it Let $D$ and $D'$ be domains in $\C$,
$D$ locally connected on $\partial D$ and $\partial D'$ weakly
flat. If $f$ is a homeomorphism with finite distortion of $D$ onto
$D'$ with $K_f\in L^{1}(D)$, then $f^{-1}$ has an extension by
continuity in $\lC$ to $\overline{D'}$.} \end{theorem}

{\it Proof.} By the Fubini theorem, the set
\begin{equation}\label{eqKPR6.2a}E\ =\ \{r\in(0,d):K_f|_{D(r)}\in L^{1}(D(r))\}
\end{equation} has a positive linear measure because $K_f\in L^{1}(D)$.

\begin{remark}\label{rmkKPR6.1a} It is clear from the proof that it is even
sufficient to assume in Theorem \ref{thKPR8.2} that $K_f$ is
integrable only in a neighborhood of $\partial D$. \end{remark}

Moreover, in view of Theorem \ref{thKPR3.1} we obtain by Theorem 9.2
in \cite{KR$_2$} or Theorem 9.7 in \cite{MRSY} the following
conclusion.

\begin{theorem}\label{thKPR8.3} {\it Let $D$ and $D'$ be domains in $\C$,
$D$ locally connected on $\partial D$ and $\partial D'$ weakly flat,
and let $f:D\to D'$ be a homeomorphism with finite distortion such
that the condition
\begin{equation}\label{eq8.10.15} \int\limits_{0}^{\delta(z_0)}
\frac{dr}{||\,K_f||_{1}(z_0,r)}\ =\ \infty\end{equation} holds for
all $z_0\in\partial D$ with some $\delta(z_0)\in(0,d(z_0))$ where
$d(z_0)\ =\ \sup\limits_{z\in D}\,|\,z-z_0|$
 and
\begin{equation}\label{eq8.10.17}
||\,K_f||_{1}(z_0,r)=\int\limits_{D(z_0,r)}K_f\ ds\end{equation} is
the $L_{1}$-norm of $K_f$ over $D(z_0,r)=\{z\in D:|z-z_0|=r\}=D\cap
S(z_0,r)$. Then there is an extension of $f^{-1}$ by continuity in
$\lC$ to $\overline{D'}$.} \end{theorem}

\cc
\section{On homeomorphic extension to the boundary}

Combining Lemma \ref{lem5.1} and Theorem \ref{thKPR8.3}, we obtain
the following statements.

\begin{theorem}\label{thKPR9.1} {\it Let $D$ and $D'$ be bounded domains in
$\C$ and let $f:D\to D'$ be a homeomorphism with finite distortion
in $D$. Suppose that the domain $D$ is locally connected on
$\partial D$ and that the domain $D'$ has a weakly flat boundary. If
\begin{equation}\label{eq8.11.2}\int\limits_{0}^{\delta(z_0)}
\frac{dr}{||\,K_f||_{1}(z_0,r)}\ =\ \infty\qquad\forall\
z_o\in\partial D\end{equation} for some $\delta(z_0)\in(0,d(z_0))$
where $d(z_0)\ =\ \sup\limits_{z\in D}\,|\,z-z_0|$ and
\begin{equation}\label{eq8.11.4}
||\,K_f||_{1}(z_0,r)=\int\limits_{D\cap S(z_0,r)}K_f\ ds\
,\end{equation} then $f$ has a homeomorphic extension to
$\overline{D}$.} \end{theorem}

In particular, as a consequence of Theorem \ref{thKPR9.1} we
obtain the following generalization of the well-known
Gehring-Martio theorem on a homeomorphic extension to the boundary
of quasiconformal mappings between QED domains, see \cite{GM}.

\begin{corollary}\label{thKPR9.2} {\it Let $D$ and $D'$ be bounded
domains with weakly flat boundaries in $\C$ and let $f:D\to D'$ be a
homeomorphism with finite distortion in $D$. If the condition
(\ref{eq8.11.2}) holds at every point $z_0\in\partial D$, then $f$
has a homeomorphic extension to $\overline{D}$.} \end{corollary}

By Theorem \ref{thKPR3.1} we have also the following, see Theorem
10.3 in \cite{KR$_2$} or Theorem 9.10 in \cite{MRSY}.

\begin{theorem}\label{thKPR9.3} {\it Let $D$ be a bounded domain in $\C$, $X\subset D$,
and $f:D\setminus\{X\}\to\lC$ a homeomorphism with finite
distortion. Suppose that $X$ and $C(X,f)$ are {\rm NED} sets. If
the condition (\ref{eq8.11.2}) holds at every point $z_0\in X$ for
$\delta(z_0)<\mathrm{dist}(z_0,\partial D)$ where
\begin{equation}\label{eq8.11.7} ||\,K_f||_{1}(z_0,r)=\int\limits_{|\,z-z_0|=r}K_f(z)\  |dz|\ ,\end{equation}
then $f$ has a homeomorphic extension to $D$.} \end{theorem}

\begin{remark}\label{rmkKPR9.1} In particular, the conclusion of Theorem \ref{thKPR9.3}
is valid if $X$ is a closed set with \end{remark}
\begin{equation}\label{eq8.11.13} H^{1}(X)\ =\ 0\ =\
H^{1}(C(X,f)).\end{equation}

\cc
\section{On some integral conditions}

Recall theorems on interconnections between some integral conditions
from \cite{RSY} and \cite{RSY$_1$}.

For every non-decreasing function $\Phi:[0,\infty]\to[0,\infty]$,
the {\bf inverse function} $\Phi^{-1}:[0,\infty]\to[0,\infty]$ can
be well defined by setting \begin{equation}\label{eqRSY2.1}
\Phi^{-1}(\tau)\ = \inf\limits_{\Phi(t)\geqslant\tau}t\,.
\end{equation} Here $\inf$ equal to $\infty$ if the set of
$t\in[0,\infty]$ such that $\Phi(t)\geqslant\tau$ is empty. Note
that the function $\Phi^{-1}$ is non-decreasing, too.

Further, the integral in (\ref{eq333F}) is understood as the
Lebesgue--Stieltjes integral and the integrals in (\ref{eq333Y}) and
(\ref{eq333B})--(\ref{eq333A}) as the ordinary Lebesgue integrals.
In (\ref{eq333Y}) and (\ref{eq333F}) we complete the definition of
integrals by $\infty$ if $\Phi(t)=\infty$, correspondingly,
$H(t)=\infty$, for all $t\geqslant T\in[0,\infty)$.

\begin{theorem}\label{thKPR3.17} {\it Let $\Phi:[0,\infty ]\to[0,\infty]$ be a
non-decreasing function and set \begin{equation}\label{eq333E} H(t)\
=\ \log \Phi(t)\,.\end{equation}

Then the equality \begin{equation}\label{eq333Y}
\int\limits_{\Delta}^{\infty}H'(t)\,\frac{dt}{t}\ =\
\infty\end{equation} implies the equality
\begin{equation}\label{eq333F} \int\limits_{\Delta}^{\infty}
\frac{dH(t)}{t}\ =\ \infty\end{equation} and (\ref{eq333F}) is
equivalent to \begin{equation}\label{eq333B}
\int\limits_{\Delta}^{\infty}H(t)\,\frac{dt}{t^2}\ =\ \infty
\end{equation} for some $\Delta>0$, and (\ref{eq333B}) is equivalent to
every of the equalities: \begin{equation}\label{eq333C}
\int\limits_{0}^{\delta}H\left(\frac{1}{t}\right)\,{dt}\ =\ \infty
\end{equation} for some $\delta>0$, \begin{equation}\label{eq333D}
\int\limits_{\Delta_*}^{\infty} \frac{d\eta}{H^{-1}(\eta)}\ =\
\infty \end{equation} for some $\Delta_*>H(+0)$,
\begin{equation}\label{eq333A} \int\limits_{\delta_*}^{\infty}
\frac{d\tau}{\tau\Phi^{-1}(\tau)}\ =\ \infty \end{equation} for some
$\delta_*>\Phi(+0)$.

Moreover, (\ref{eq333Y}) is equivalent  to (\ref{eq333F}) and hence
(\ref{eq333Y})--(\ref{eq333A}) are equivalent each to other  if
$\Phi$ is in addition absolutely continuous. In particular, all the
conditions (\ref{eq333Y})--(\ref{eq333A}) are equivalent if $\Phi$
is convex and non-decreasing.} \end{theorem}

\begin{remark}\label{rmk1} It is necessary to give one more explanation.
From the right hand sides in the conditions
(\ref{eq333Y})--(\ref{eq333A}) we have in mind $+\infty$. If
$\Phi(t)=0$ for $t\in[0,t_*]$, then $H(t)=-\infty$ for $t\in[0,t_*]$
and we complete the definition $H'(t)=0$ for $t\in[0,t_*]$. Note,
the conditions (\ref{eq333F}) and (\ref{eq333B}) exclude that $t_*$
belongs to the interval of integrability because in the contrary
case the left hand sides in (\ref{eq333F}) and (\ref{eq333B}) are
either equal to $-\infty$ or indeterminate. Hence we may assume in
(\ref{eq333Y})--(\ref{eq333C}) that $\Delta>t_0$ where $t_0\colon
=\sup\limits_{\Phi(t)=0}t$, $t_0=0$ if $\Phi(0)>0$, and
$\delta<1/t_0$, correspondingly. \end{remark}

\begin{theorem}\label{th5.555} {\it Let $Q:{\Bbb D}\to[0,\infty]$ be a measurable
function such that \begin{equation}\label{eq5.555} \int\limits_{\Bbb
D}\Phi(Q(z))\,dxdy\ <\ \infty\end{equation} where
$\Phi:[0,\infty]\to[0,\infty]$ is a non-decreasing convex function
such that \begin{equation}\label{eq3.333a}
\int\limits_{\delta_0}^{\infty}\ \frac{d\tau}{\tau\Phi^{-1}(\tau)}\
=\ \infty\end{equation} for some $\delta_0>\Phi(0)$. Then
\begin{equation}\label{eq3.333A} \int\limits_{0}^{1}\
\frac{dr}{rq(r)}\ =\ \infty \end{equation} where $q(r)$ is the
average of the function $Q(z)$ over the circle $|z|=r$.}
\end{theorem}

Here ${\Bbb D}$ denotes the unit disk in ${\Bbb C}$. Combining
Theorems \ref{thKPR3.17} and \ref{th5.555} we obtain also the
following.

\begin{corollary}\label{cor555} {\it If $\Phi:[0,\infty]\to[0,\infty]$ is a
non-decreasing convex function and $Q:{\Bbb D}\to[0,\infty]$
satisfies (\ref{eq5.555}), then every of the conditions
(\ref{eq333Y})--(\ref{eq333A}) implies (\ref{eq3.333A}).}
\end{corollary}

\cc
\section{On the mappings quasiconformal in the mean}

Integral conditions of the type
\begin{equation}\label{eq1.1KR}\int\limits_{D}\Phi(K(x))\,dm(x)<\infty\end{equation}
are often applied in the mapping theory, see e.g. \cite{Ah},
\cite{Bi}, \cite{Gol}, \cite{Kr}--\cite{Ku}, \cite{Pes},
\cite{Rya} and \cite{UV}.

Combining Theorem \ref{th5.555} with Lemma \ref{lem5.1} and Theorem
\ref{thKPR9.1}, we come to the following statement.\newpage

\begin{theorem}{}\label{thKR4.1} {\it Let $D$ and $D'$ be bounded domains in $\C$
such that $D$ is locally connected at $\partial D$ and $D'$ has a
weakly flat (strongly accessible) boundary. Suppose that $f:D\to
D'$ is a homeomorphism with finite distortion and
\begin{equation}\label{eqKR4.1}\int\limits_{D}\Phi(K_{f}(z))\,dm(z)<\infty\end{equation}
for a convex non-decreasing function
$\Phi:[0,\infty]\to[0,\infty]$. If
\begin{equation}\label{eqKR4.2}\int\limits_{\delta_0}^{\infty}\frac{d\tau}{\tau\Phi^{-1}(\tau)}=
\infty\end{equation} for some $\delta_0>\Phi(0)$, then $f$ has a
homeomorphic (continuous) extension $\overline{f}$ to $\overline{D}$
that maps $\overline{D}$ onto (into) $\overline{D'}$.} \end{theorem}

\begin{remark}\label{rmKR4.1} In particular, the conclusion on homeomorphic extension
is valid for domains $D$ and $D'$ with smooth boundaries and for
convex domains. Note also that by Theorem \ref{thKPR3.17} the
condition (\ref{eqKR4.2}) can be replaced by each of the conditions
(\ref{eq333Y}) -- (\ref{eq333D}).  The example in \cite{KR$_3$}
shows that each of the given conditions are not only sufficient but
also necessary for continuous extension of $f$ to the boundary.
\end{remark}

\medskip
\noindent Kovtonyk D., Petkov I. and Ryazanov V., \\
Institute of Applied Mathematics and Mechanics,\\
National Academy of Sciences of Ukraine, \\
74 Roze Luxemburg str., 83114 Donetsk, UKRAINE \\
Phone: +38 -- (062) -- 3110145, Fax: +38 -- (062) -- 3110285 \\
denis$\underline{\ \ }$\,kovtonyuk@bk.ru, igorpetkov@list.ru,
vlryazanov1@rambler.ru

\end{document}